# RESEARCH ANNOUNCEMENT

# CHAOS IN THE LORENZ EQUATIONS: A COMPUTER-ASSISTED PROOF

KONSTANTIN MISCHAIKOW AND MARIAN MROZEK

ABSTRACT. A new technique for obtaining rigorous results concerning the global dynamics of nonlinear systems is described. The technique combines abstract existence results based on the Conley index theory with computer-assisted computations. As an application of these methods it is proven that for an explicit parameter value the Lorenz equations exhibit chaotic dynamics.

## INTRODUCTION

The purpose of this note is to briefly describe a new technique for obtaining rigorous results concerning the global dynamics of nonlinear systems. The technique combines abstract existence results based on topological invariants (the Conley index) with finite, computer-assisted computations necessary to verify the assumptions of the theorems in a concrete example. There are at least three new aspects of this technique. It applies to concrete differential equations; it can provide a relatively strong description of global dynamics (in terms of semiconjugacies); and the necessary computer-assisted computations are small enough to be performed on currently available computers.

To focus the presentation of the ideas, an outline of a proof that the Lorenz equations,

(1)
$$\begin{aligned}
\dot{x} &= s(y-x), \\
\dot{y} &= Rx - y - xz, \\
\dot{z} &= xy - qz,
\end{aligned}$$

contain chaotic dynamics for a prescribed open set of parameter values will be presented.

Let $f : \mathbf{R}^n \to \mathbf{R}^n$ be a homeomorphism. For $N \subset \mathbf{R}^n$, the *maximal invariant set* of $N$ is defined by $\operatorname{Inv}(N, f) = \{x \in N \mid f^n(x) \in N \ \forall \ n \in \mathbf{Z}\}$.

Received by the editors April 10, 1993, and, in revised form, February 15, 1994.
1991 *Mathematics Subject Classification*. Primary 58F13, 54820, 65L99.
The first author's research was funded in part by NSF Grants DMS-9101412 and DMS-9302970.
The second author is on leave from the Computer Science Department, Jagiellonian University, Kraków, Poland.

**Theorem 1.** *Let*
$$P := \{(x, y, z) \mid z = 53\}.$$
*For all parameter values in a sufficiently small neighborhood of* $(s, R, q) = (45, 54, 10)$, *there exists a Poincare section* $N \subset P$ *such that the Poincare map* $g$ *induced by* (1) *is Lipschitz and well defined. Furthermore, there exists a* $d \in \mathbf{N}$ *and a continuous surjection* $\rho : \mathrm{Inv}(N, g) \to \Sigma_2$ *such that*
$$\rho \circ g^d = \sigma \circ \rho$$
*where* $\sigma : \Sigma_2 \to \Sigma_2$ *is the full shift dynamics on two symbols.*

The choice of parameters in Theorem 1 was dictated mainly by an attempt to minimize the necessary computations. Preliminary estimates show that analogous computations for the classical choice of parameters $(s, R, q) = (10, 28, 8/3)$, though substantially more complex, are also in the range of currently available computers (computations in progress, see [4]). However, it must be emphasized that Theorem 1 shows only the existence of an unstable invariant set which maps onto a horseshoe. Though the set may lie within a strange attractor, we do not yet have sufficiently strong abstract results to prove that the whole attractor is chaotic.

Let us mention some other attempts to reduce the question of chaos in the Lorenz equations to a finite computation: an $\Omega$-explosion approach [7] and a shooting method approach [1]. Contrary to the presented method, in both cases the authors have not yet performed the actual computations.[1] Moreover, as mentioned before, Theorem 1 provides a description of the dynamics in a neighborhood of an explicitly presented parameter value, i.e. $(s, R, q) = (45, 54, 10)$. This should be contrasted with the other methods where the conclusion is that explicit dynamics occur for some unknown parameter value within a specified range (shooting method) or for a sequence of values tending to a given parameter value (the $\Omega$-explosion method).

The proof of Theorem 1 has five distinct components:

1. algebraic invariants based on the Conley index theory which guarantee the structure of the global dynamics, in this case the semi-conjugacy to the full two shift;
2. an extension of these invariants to multivalued maps;
3. a theory of finite representable multivalued maps which, when combined with the above-mentioned algebraic invariants, serves to bridge the gap between the continuous dynamics (in this case the Lorenz equations) and the finite dynamics of the computer;
4. the numerical computations of the finite multivalued map of interest;
5. the combinatorial computations of the Conley index for the multivalued map.

The rough idea of the general scheme of the proof is as follows. Choose a potential isolating neighborhood $N$ for the Poincare map $g$. Select a finite representable multivalued map $\mathcal{G}$ (definitions follow) such that $\mathcal{G}$ is an *extension* of $g$, i.e. $g(x) \in \mathcal{G}(x)$ for all $x \in \mathrm{dom}\,\mathcal{G}(x)$. Perform a computer calculation both to determine $\mathcal{G}$ and to check whether $N$ is an isolating neighborhood for $\mathcal{G}$. Theorem 3 guarantees that whenever $N$ is an isolating neighborhood for $\mathcal{G}$, it is also an isolating neighborhood for $g$ and the Conley indexes coincide. If $N$ is an isolating neighborhood, then the Conley index of $N$ under $\mathcal{G}$ is computed (this is again a finite computation and can be done by computer). This then determines the Conley index for the set $N$ under the original Poincare map $g$ and allows one to verify the assumptions of Theorem 2, which yields Theorem 1.

---

[1]Since completion of the original version of this paper, numerical calculations based on modifications of the ideas of [1] have been completed [8].

If $N$ is not an isolating neighborhood of $\mathcal{G}$, then one is free to choose a new multivalued map which is a better approximation of $g$ and to repeat the computations. Theorem 4 implies that if $N$ is an isolating neighborhood for $g$, then $N$ is an isolating neighborhood for a sufficiently refined choice of $\mathcal{G}$. In other words, if $N$ is an isolating neighborhood under $g$, then given sufficient computing power, this scheme will provide a proof of this fact and, in addition, the values of the Conley index of $N$ under $g$.

In the sections below each of the five components of the proof will be discussed. Most of the theoretical results hold true in greater generality; however, for the sake of exposition the simplest setting has been chosen.

## 1. Isolating neighborhoods and chaos

Again, let $f : \mathbf{R}^n \to \mathbf{R}^n$ be a homeomorphism. A compact set $N$ is called an *isolating neighborhood* if $\text{Inv}(N, f) \subset \text{int } N$. To determine the Conley index of an isolated invariant set for a map, one begins with an *index pair* $(N, L)$. The details of the definitions which follow can be found in [5]. For the purposes of this paper it is sufficient to know that $N$ is an isolating neighborhood, $L$ is its exit set, and the map $f$ induces a homomorphism called the *index map* on the cohomology of the index pair, i.e. $I_f^* : H^*(N, L) \to H^*(N, L)$. The cohomological Conley index of an isolating neighborhood under $f$ is given by

$$\text{Con}^*(\text{Inv}(N, f)) = \Big(CH^*(\text{Inv}(N, f)), \chi^*(\text{Inv}(N, f))\Big),$$

where $CH^*(\text{Inv}(N, f))$ is the graded module obtained by quotienting $H^*(N, L)$ by the generalized kernel of $I_f^*$ and $\chi^*(\text{Inv}(N, f))$ is the induced graded module automorphism on $CH^*(\text{Inv}(N, f))$.

Let $N = N_0 \cup N_1$ be an isolating neighborhood under $f$ where $N_0$ and $N_1$ are disjoint compact sets. For $k, l = 0, 1$, let $N_{kl} = N_k \cap f(N_l)$. Let $S_k := \text{Inv}(N_{kk}, f)$ and $S_{lk} := \text{Inv}(N_{kk} \cup N_{kl} \cup N_{ll}, f)$. The following result is a special case of [3, Theorem 2.3].

**Theorem 2.** *Assume that*

$$\text{Con}^n(S_k) = \begin{cases} (\mathbf{Q}, \text{id}) & \text{if } n = 1, \\ 0 & \text{otherwise} \end{cases}$$

*and that $\chi^*(S_{lk})$ is not conjugate to $\chi^*(S_k) \oplus \chi^*(S_l)$. Then there exists $d \in \mathbf{N}$ and a continuous surjection $\rho : \text{Inv}(N, f) \to \Sigma_2$ such that*

$$\rho \circ f^d = \sigma \circ \rho$$

*where $\sigma : \Sigma_2 \to \Sigma_2$ is the full shift dynamics on two symbols.*

## 2. Index theory for multivalued maps

Recall that a *multivalued map* from $\mathbf{R}^n$ to itself is a function $\mathcal{F} : \mathbf{R}^n \to \mathcal{P}(\mathbf{R}^n) \setminus \{\varnothing\}$ from $\mathbf{R}^n$ to the power set of $\mathbf{R}^n$. A continuous function $f : \mathbf{R}^n \to \mathbf{R}^n$ is a *selector* for $\mathcal{F}$ if $f(x) \in \mathcal{F}(x)$ for all $x \in \mathbf{R}^n$. In the other direction, a multivalued function $\mathcal{F}$ is an *extension* of a continuous function $f : \mathbf{R}^n \to \mathbf{R}^n$ if $f(x) \in \mathcal{F}(x)$ for all $x \in \mathbf{R}^n$. Let $\mathcal{F}$ be a multivalued function on $\mathbf{R}^n$. For $A \subset \mathbf{R}^n$ let $\mathcal{F}(A) := \bigcup_{x \in A} \mathcal{F}(x)$ and define, recursively, $\mathcal{F}^{n+1}(A) := \mathcal{F}(\mathcal{F}^n(A))$. In this way $\mathcal{F}$ defines a multivalued discrete semidynamical system on $\mathbf{R}^n$. A Conley index theory for such systems has been developed in [2], and the basic ideas are as follows. Given $B \subset \mathbf{R}^n$, its *inverse image* is $\mathcal{F}^{-1}(B) := \{x \in \mathbf{R}^n \mid \mathcal{F}(x) \subset B\}$ and its *weak inverse image* is $\mathcal{F}^{*-1}(B) := \{x \in \mathbf{R}^n \mid \mathcal{F}(x) \cap B \neq \varnothing\}$. Given $N \subset \mathbf{R}^n$, the invariant set of $N$ is given by

$$\mathrm{Inv}(N, \mathcal{F}) := \{x \in N \mid \exists \gamma_x : \mathbf{Z} \to N \text{ such that } \gamma_x(0) = 0, \text{ and}$$

$$\gamma_x(n+1) \in \mathcal{F}(\gamma(n))\}.$$

The *diameter* of $\mathcal{F}$ over $N$ is the number

$$\mathrm{diam}_N \mathcal{F} := \sup_{x \in N} \{\|z - y\| \mid z, y \in \mathcal{F}(x)\}.$$

The diameter of $\mathcal{F}$ over its domain will also be called the *size* of $\mathcal{F}$. $N$ is called an *isolating neighborhood* under $\mathcal{F}$ if

$$B(\mathrm{Inv}(N, \mathcal{F}), \mathrm{diam}_N \mathcal{F}) \subset \mathrm{int}\,(N).$$

In order to define the Conley index for $\mathcal{F}$ in $N$, some conditions (admissibility) must be met (see [2]). For the sake of this paper it is enough to mention that convex-valued maps with continuous selectors are admissible.

**Theorem 3.** *Let $\mathcal{F}$ be an admissible multivalued map and $f$ a selector of $\mathcal{F}$. If $N$ is an isolating neighborhood for $\mathcal{F}$, then it is an isolating neighborhood for $f$. Furthermore,*

$$\mathrm{Con}^*\bigl(\mathrm{Inv}(N, f)\bigr) \approx \mathrm{Con}^*\bigl(\mathrm{Inv}(N, \mathcal{F})\bigr).$$

The importance of this theorem is that it implies that if the Poincare map $g$ given by the Lorenz equations is replaced by an extension $\mathcal{G}$ (i.e. a multivalued map such that $g(x) \in \mathcal{G}(x)$) such that $N$ is an isolating neighborhood for $\mathcal{G}$, then any index information obtained for $\mathcal{G}$ is valid for $g$.

The next theorem is a typical convergence theorem and states that any index information of $g$ can be determined by a sufficiently small extension $\mathcal{G}$.

**Theorem 4.** *Let $N$ be an isolating neighborhood for $f : \mathbf{R}^n \to \mathbf{R}^n$, a Lipschitz continuous function. Let $\{\mathcal{F}_n\}$ be a family of extensions of $f$ such that $\mathcal{F}_n \to f$. Then for $n$ sufficiently large, $N$ is an isolating neighborhood for $\mathcal{F}$.*

For the purposes of computation it is convenient to use the following special form of an isolating neighborhood. $N$ is an *isolating block* for $\mathcal{F}$ if

(2) $$B\bigl(\mathcal{F}^{*-1}(N) \cap N \cap \mathcal{F}(N), \mathrm{diam}_N \mathcal{F}\bigr) \subset \mathrm{int}\,(N).$$

Notice that contrary to the notion of an isolating neighborhood, it uses only a finite number of iterates of $\mathcal{F}$ (one forward and one backward).

## 3. Finite representable multivalued maps

For the proof it is necessary to use a computer (which can manipulate only a finite set of objects) to obtain cohomological information generated by the continuous map $g$. Since $g$ is Lipschitz, the simplicial approximation theorem guarantees that in principle this can be done but the computations must be performed in such a way that they apply not only to the simplicial approximation but also to all nearby maps, including the original map $g$. Hence it is more natural and technically simpler to employ finite multivalued maps. The ideas behind this solution will now be explained (see [6] for details).

On a computer one can work with infinite sets of mathematical objects (numbers, vectors, functions, relations, etc.) only by choosing a finite subset of objects and coding its elements with natural numbers. The elements of the selected finite subset of objects are referred to as representable objects (relative to the selected coding). In the proof of Theorem 1 the set of representable real numbers is given by the standard floating point representable coordinates. The representable sets in $\mathbf{R}^3$ were chosen by selecting a compact set $M \subset \mathbf{R}^3$ within which the dynamics of interest occur, a set $D$ of representable vectors, and a representable number $\eta$, in such a way that the set $M$ is covered by the collection of balls $B(d,\eta)$, $d \in D$. Since the particular metric chosen influences only the efficiency of the algorithm, for the proof the sup norm was chosen. Thus, the balls are in fact cubes. A *representable* set is any union of a subcollection of these cubes. Let $M_0 \subset M$. A *representable multivalued map* on $M_0$ is a multivalued map $\mathcal{F} : M_0 \to \mathcal{P}(M)$ such that the set $\{\mathcal{F}(x) \mid x \in M_0\}$ is a finite collection of representable sets.

Representable multivalued maps allow for the passage from continuous dynamics of the Lorenz equations to the discrete dynamics of the computer. Furthermore, they are more readily attainable than simplicial approximations, as can be seen from the following description. Let $f : M_0 \to M$ be a Lipschitz continuous function with Lipschitz constant $L$. Assume that for every representable vector $d$ one can compute a representable vector $f_0(d)$ such that $\|f(d) - f_0(d)\| < \delta$, where $\delta$ is a given representable number. Defining $\mathcal{F}(x)$ as the smallest convex representable set which contains $B(f_0(x), \delta + L\eta)$, one obtains a multivalued map $\mathcal{F} : M_0 \to M$ such that $f(x) \in \mathcal{F}(d)$ for all $x \in B(d, \eta)$. Then

$$\mathcal{F}^u(x) := \bigcup_{\|d-x\| \leq \eta} \mathcal{F}(d), \text{ and } \mathcal{F}^l(x) := \bigcap_{\|d-x\| \leq \eta} \mathcal{F}(d)$$

are easily seen to be finite, representable, and respectively upper and lower semi-continuous extensions of $f$. Furthermore, $\mathcal{F}^l$ is convex valued.

Letting $\delta \to 0$ results in the following proposition.

**Proposition 5.** *Let $M_0 \subset \mathbf{R}^n$ be a compact set, and let $f : M_0 \to \mathbf{R}^n$ be a Lipschitz continuous function. Then there exists a sequence $\mathcal{F}_n$ of finite representable extensions of $f$ such that $F_n \to f$ as $n \to \infty$.*

## 4. Numerical computations

Simple and standard numerical integration of the Lorenz equations at the parameter values (45, 54, 10) leads to a return map which strongly suggests the existence of horseshoe dynamics. In particular, one can find two rectangles $R_0$ and $R_1$ in the plane $P$ which, under the numerical integration, appear to cross themselves and each other transversally. Thus, one expects the existence of an invariant set which is conjugate to the full shift dynamics on two symbols.

As was suggested in the introduction, the proof of the semi-conjugacy described in Theorem 1 was obtained by computing a representable multivalued admissible extension $\mathcal{G}$ of the Poincare map $g$ such that $N := N_0 \cup N_1$ is an isolating neighborhood under $\mathcal{G}$. The sets $N_0$ and $N_1$ were taken to be carefully chosen rectangles within $R_0$ and $R_1$ in order to minimize computation. Ideally the multi-valued extension $\mathcal{G}$ would have been obtained by integrating the Lorenz equations numerically from the center of each cube in the grid of representable sets and incorporating all errors into the size of the assigned value. Unfortunately the growth of errors is exponential with time, which caused technical difficulties that needed to be overcome.

There are four sources of error:
   (i) the approximate arithmetic of the machine;
   (ii) the numerical procedure used to integrate the equations;
   (iii) extending the value from the center of a cube to the whole cube by means of Lipschitz constant estimates, as was explained in the previous section;
   (iv) estimating the point of intersection of the trajectory with the cross-section from two consecutive steps of the numerical method.

In every case rigorous error bounds can be obtained (see [4] for details). It should be mentioned here that double precision arithmetic and the standard fourth-order Runge-Kutta method with step size $100/2^{20}$ were used. The error estimates for (iii) were based on the Gronwall inequality, the local Lipschitz constants, and logarithmic norms. The second-order Taylor expansion of the solution was used to estimate errors in (iv). This resulted in the growth factor (ratio of the size of a value to the size of the grid) for equation (1) and in the function $\mathcal{G}$ described above being approximately $10^6$. The number of initial grid points needed to fulfill condition (2) by the associated multivalued map $\mathcal{G}$ is proportional to the square of the growth factor.

With the computing power available a significantly smaller growth factor was needed. To achieve this, twenty-three intermediate cross-sections, labelled $\Xi_k$, $k = -1, \ldots, 22$, were introduced, with the original Poincare section appearing as $P = \Xi_{-1} = \Xi_{22}$. Using the technique described above, twenty-two finite representable multivalued maps $\mathcal{G}_i : \Xi_{i-1} \to \Xi_i$ were obtained. Observe that each multivalued map is an extension of the flow defined map $g_i : \Xi_{i-1} \to \Xi_i$ determined by the Lorenz equations. $\mathcal{G}$ was taken to be the composition $\mathcal{G}_{21} \circ \mathcal{G}_{20} \circ \ldots \circ \mathcal{G}_0$.

The calculation was performed beginning with approximately 700,000 cubes covering $N_1$. (It was not necessary to compute $\mathcal{G}$ on $N_0$ because of the symmetry present in the Lorenz equations.) The growth factor in this case was approximately 30, which resulted in the diameter of $\mathcal{G}(q)$ being less than 0.044 for all $q \in Q$. Thus to determine whether $N$ is an isolating block under $\mathcal{G}$, one need only check that
$$B\big(\mathcal{G}^{*-1}(N) \cap N \cap \mathcal{G}(N), 0.044\big) \subset \text{int } N.$$
This is indeed the case. Therefore $N$ is an isolating neighborhood of $\mathcal{G}$, and hence by Theorem 3 $N$ is an isolating neighborhood of $g$.

## 5. The Conley index of $\mathcal{G}$

Recall that to compute the Conley index of an isolated invariant set for a map, one first needs to find an index pair $(N, L)$ and then to determine the index map, i.e.
$$I_{\mathcal{G}}^* : H^*(N, L) \to H^*(N, L).$$
For this particular index pair these objects can be determined essentially by inspection, since the only cohomology group of interest is $H^1(N, L)$. However, for more general problems it should be observed that $D$ can be viewed as the set of vertices of a simplicial decomposition of $P$, and hence $\mathcal{G}$ is a multivalued simplicial map. Therefore, determining $I_\mathcal{G}$ is in fact a finite computation and hence in principle can be done by the computer. From $I_\mathcal{G}^*$ one computes the Conley index $\mathrm{Con}^*(\mathrm{Inv}(N, \mathcal{G}))$ and hence by Theorem 3 the index of interest $\mathrm{Con}^*(\mathrm{Inv}(N, g))$. For this problem the index computations lead to the algebraic hypotheses of Theorem 2, and hence Theorem 1 follows.

## Acknowledgments

Both authors gratefully acknowledge Luca Dieci's advice regarding the numerical computations necessary in this work. It is a pleasure for the second author to thank Colin Sparrow for several conversations held many years ago, which turned out to be crucial in the choice of the Lorenz equations as an illustration of the presented method.

School of Mathematics, Georgia Institute of Technology, Atlanta, Georgia 30332

*E-mail address*, M. Mrozek: `mrozek@ii.uj.edu.pl`

*E-mail address*, K. Mischaikov: `mischaik@math.gatech.edu`